\documentclass[10pt,twoside,reqno,psamsfonts]{amsart}

\usepackage[OT1]{fontenc}
\usepackage{type1cm}
\usepackage{enumerate}
\usepackage{amssymb}
\usepackage[all]{xy}
\usepackage{mathrsfs}
\usepackage[dvips]{graphicx}
\usepackage{psfrag}          
\usepackage{geometry}        
\usepackage{version}         

\numberwithin{equation}{section}

\theoremstyle{plain}
\newtheorem{theorem}{Theorem}[section]

\newtheorem{proposition}[theorem]{Proposition}
\newtheorem{lemma}[theorem]{Lemma}

\theoremstyle{remark}
\newtheorem{remark}[theorem]{Remark}

\newtheorem*{ack}{Acknowledgements}

\theoremstyle{definition}
\newtheorem{definition}[theorem]{Definition}

\newenvironment{note}{\begin{quote}\sf}{\end{quote}}

\newcommand{\sE}{\mathscr{E}}

\newcommand{\R}{\mathbb{R}}

\newcommand{\N}{\mathbb{N}}

\DeclareMathOperator{\diam}{diam}

\DeclareMathOperator{\spt}{spt}

\DeclareMathOperator{\Geo}{Geo}
\DeclareMathOperator{\GeoOpt}{GeoOpt}

\def\XXint#1#2#3{{\setbox0=\hbox{$#1{#2#3}{\int}$ }
\vcenter{\hbox{$#2#3$ }}\kern-.6\wd0}}

\begin{document}

\title[Improved geodesics for $CD^*(K,N)$]
{Improved geodesics for the reduced curvature-dimension condition in branching metric spaces}

\author{Tapio Rajala}
\address{Scuola Normale Superiore\\
Piazza dei Cavalieri 7\\
I-56127 Pisa\\ Italy}
\email{tapio.rajala@sns.it}

\thanks{The author acknowledges the support of the European Project ERC AdG *GeMeThNES* and the Academy of Finland project no. 137528.}
\date{\today}


\begin{abstract}
 In this note we show that in metric measure spaces satisfying the reduced curvature-dimension condition $CD^*(K,N)$ we always 
 have geodesics in the Wasserstein space of probability measures that satisfy the critical convexity inequality of $CD^*(K,N)$
 also for intermediate times and in addition the measures along these geodesics have an upper-bound on their densities. 
 This upper-bound depends on the bounds for the densities of the end-point measures, the lower-bound $K$ for the Ricci-curvature,
 the upper-bound $N$ for the dimension, and on the diameter of the union of the supports of the end-point measures.
\end{abstract}


\maketitle

\section{Introduction}

Different definitions for lower Ricci curvature bounds in metric measure spaces have been studied extensively in the recent years.
A common goal for these definitions is to extend the Riemannian definition to a subset of metric measure spaces which is closed
under the measured Gromov-Hausdorff convergence. 
Among the most celebrated definitions are the ones given by Sturm \cite{S2006I, S2006II} and by Lott and Villani \cite{LV2009} using optimal mass
transportation.
In both of these treatments the critical convexity property for the definitions involving an upper-bound on the dimension was
given using the R\'enyi entropy $\sE_N \colon \mathcal{P}(X) \to [-\infty,0]$ which is defined as
\[
 \sE_N(\mu) = -\int_X \rho^{1-1/N} \,dm,
\]
where $\rho$ is the density of $\mu$ with respect to the measure $m$ which is the underlying measure in our
metric measure space $(X,d,m)$. In this note we will always assume the metric space $(X,d)$ to be complete and separable, 
and the measure $m$ to be locally finite.

A definition by Sturm \cite{S2006II} is that a metric measure space $(X,d,m)$ is a $CD(K,N)$-space,
with the interpretation that it has $N$-Ricci curvature bounded below by $K$, if for any two measures 
$\mu_0, \mu_1 \in \mathcal{P}(X)$ with $W_2(\mu_0,\mu_1)<\infty$ there exists $\pi \in \GeoOpt(\mu_0,\mu_1)$ so that
along the Wasserstein geodesic $\mu_t = (e_t)_\#\pi$ for every $t \in [0,1]$ and $N' \ge N$ we have
\begin{equation}\label{eq:CD-def}
 \sE_{N'}(\mu_t) \le - \iint_{X\times X}\left[\tau_{K,N'}^{(1-t)}(d(x_0,x_1))\rho_0^{-1/N'}(x_0) 
                 + \tau_{K,N'}^{(t)}(d(x_0,x_1))\rho_1^{-1/N'}(x_1)\right]\,dq(x_0,x_1),
\end{equation}
where we have written $\mu_0 = \rho_0m + \mu_0^s$ and $\mu_1 = \rho_1m + \mu_1^s$ with $\mu_0^s \perp m$, $\mu_1^s \perp m$ and 
$q = (e_0,e_1)_\#\pi$.
We will recall some basic notation and facts about the Wasserstein space $(\mathcal{P}(X),W_2)$ in Section \ref{sec:W}.
The volume distortion coefficients $\tau_{K,N}^{(t)}(r)$ in the inequality \eqref{eq:CD-def} are defined
as 
\[
\tau_{K,N}^{(t)}(r) = t^{1/N}\sigma_{K,N-1}^{(t)}(r)^{1-1/N}, 
\]
where
\[
 \sigma_{K,N}^{(t)}(r) = \begin{cases}
                          \infty, & \text{if }Kr^2 \ge N\pi^2 \\
                          t \frac{\mathfrak{S}_{K/N}(tr)}{\mathfrak{S}_{K/N}(r)} & \text{otherwise}
                         \end{cases}
\]
and
\[
 \mathfrak{S}_k(r) = \begin{cases}
                      \frac{\sin(\sqrt{k}r)}{\sqrt{k}r}, & \text{if }k > 0 \\
                       1, & \text{if }k = 0\\
                      \frac{\sinh(\sqrt{-k}r)}{\sqrt{-k}r}, & \text{if }k < 0.
                     \end{cases}
\]

It turns out that a slightly weaker version of $CD(K,N)$ where the coefficients $\tau_{K,N}^{(t)}(r)$ are replaced by
the smaller coefficients $\sigma_{K,N}^{(t)}(r)$ has properties which are not known for the original $CD(K,N)$-spaces.
This weaker version is called the reduced curvature-dimension condition $CD^*(K,N)$.
It was introduced and studied by Bacher and Sturm in \cite{BS2010}.
Notice that in their definition Bacher and Sturm required the convexity-inequality to hold only between any two absolutely continuous measures
with bounded supports.

\begin{definition}
 A metric measure space $(X,d,m)$ satisfies the \emph{reduced curvature-dimension condition} $CD^*(K,N)$ if for any two measures 
 $\mu_0, \mu_1 \in \mathcal{P}^{ac}(X,m)$ with bounded supports there exists a measure $\pi \in \GeoOpt(\mu_0,\mu_1)$ so that
 for $\mu_t = (e_t)_\#\pi$ we have
 \begin{equation}\label{eq:CDred-def}
 \sE_{N'}(\mu_t) \le - \iint_{X\times X}\left[\sigma_{K,N'}^{(1-t)}(d(x_0,x_1))\rho_0^{-1/N'}(x_0) 
                 + \sigma_{K,N'}^{(t)}(d(x_0,x_1))\rho_1^{-1/N'}(x_1)\right]\,dq(x_0,x_1),
 \end{equation}
 for every $t \in [0,1]$ and $N' \ge N$.
\end{definition}

In \cite{BS2010, CS2011, DS2011} the $CD^*(K,N)$-spaces were studied mainly under the extra assumption that the space is nonbranching,
meaning that there are no geodesics that start as the same geodesic but then branch out.
The main results of \cite{BS2010} include the local-to-global and tensorization properties of nonbranching $CD^*(K,N)$-spaces. 
These two properties are perhaps the two most important properties which are not known for the usual $CD(K,N)$-spaces.
However, a disadvantage of the $CD^*(K,N)$ definition when
compared to $CD(K,N)$ is that it does not immediately give the sharp constants (of the Riemannian setting) for many 
geometric inequalities. Only very recently some of the inequalities, namely the Bishop-Gromov volume growth inequality
and the (global) Poincar\'e inequality were proven in their sharp form in nonbranching $CD^*(K,N)$-spaces.
This was a byproduct of the sharp measure contraction property proven by Cavalletti and Sturm in \cite{CS2011}.
It would be very interesting to know if all these results (local-to-global, tensorization, measure contraction property,
sharp Bishop-Gromov and Poincar\'e inequalities) really require the extra assumption on nonbranching. 
It is important to keep in mind that the nonbranching assumption is not stable under measured Gromov-Hausdorff 
limits, not even when it is combined with the curvature-dimension conditions mentioned in this note.

The main advantage of the nonbranching assumption is that under it the curvature-dimension condition gives local information
about the space: the curvature-inequalities hold even along individual geodesics and at intermediate times. Without the 
nonbranching assumption information along individual geodesics seems to be out of reach. However, as will be seen in this note we can get
information on the curvature for intermediate times also in branching $CD^*(K,N)$-spaces. The reason behind this is the way how the
coefficients $\sigma_{K,N}^{(t)}(r)$ behave under localization. Basically, this is also the reason why the local-to-global 
property holds in nonbranching $CD^*(K,N)$-spaces.

Without the localization of the curvature-dimension condition to geodesics in the space $(X,d)$ we are forced to work 
with geodesics in the Wasserstein space $(\mathcal{P}(X), W_2)$. The goal of this note is to help in this regard by
proving the existence of good Wasserstein geodesics in the possibly branching $CD^*(K,N)$-spaces. In \cite{R2011b} we studied
the same question for the branching $CD(K,N)$-spaces and obtained a local Poincar\'e inequality and the measure contraction
property using the existence of good geodesics along which we had a bound on the densities.
In $CD^*(K,N)$-spaces we are able to get even better geodesics.
The geodesics which we get do not only have a similar upper-bound for the densities, but they also 
satisfy the curvature-dimension inequality of $CD^*(K,N)$ for intermediate times.

Under the nonbranching assumption it is known that \eqref{eq:CDred-def} holds for all $\pi \in \GeoOpt(\mu_0,\mu_1)$.
This implies, for instance, that in nonbranching $CD^*(K,N)$-spaces for any $\pi \in \GeoOpt(\mu_0,\mu_1)$ we have
for all $0 \le s < t < r \le 1$ and $N' \ge N$ the inequality
 \begin{equation}\label{eq:CD-inter}
 \sE_{N'}(\mu_t) \le - \iint_{X\times X}\left[\sigma_{K,N'}^{(\frac{r-t}{r-s})}(d(x_s,x_r))\rho_s^{-1/N'}(x_s) 
                 + \sigma_{K,N'}^{(\frac{t-s}{r-s})}(d(x_s,x_r))\rho_r^{-1/N'}(x_r)\right]\,dq(x_s,x_r),
 \end{equation}
where $q = (e_s,e_r)_\#\pi$. A similar inequality is also valid in nonbranching $CD(K,N)$-spaces. It is not known if
this is true in $CD(K,N)$-spaces without the nonbranching assumption. In this note we will prove that the inequality \eqref{eq:CD-inter}
holds for $N' = N$ in $CD^*(K,N)$-spaces even without the nonbranching assumption.
This stronger version of the curvature-dimension inequality is extremely useful because it for example implies
\[
 \partial_t^2\sE_N(\mu_t) \ge \frac{K}{N} \int_{\Geo(X)}\rho_t^{-\frac1N}(\gamma_t)l(\gamma)^2\,d\pi(\gamma),
\]
where $l(\gamma)$ denotes the length of the geodesic $\gamma$. For negative $K$ this further implies the simpler
looking differential inequality
\[
 \partial_t^2\sE_N(\mu_t) \ge \frac{K}{N}D^2 \sE_N(\mu_t),
\]
where $D = \diam(\spt\mu_0\cup\spt \mu_1)$.

Let us state in the following Theorem \ref{thm:main} the two properties which we prove for the good geodesic 
that are constructed in this note.

\begin{theorem}\label{thm:main}
 Suppose that $(X,d,m)$ satisfies the reduced curvature-dimension condition $CD^*(K,N)$ for some $K \in \R$
 and $N \in (1,\infty)$. Then for all
 $\mu_0, \mu_1 \in \mathcal{P}^{ac}(X,m)$ with bounded densities and $D = \diam(\spt\mu_0\cup\spt \mu_1) < \infty$ there
 exists $\pi \in \GeoOpt(\mu_0,\mu_1)$ which
 \begin{enumerate}
  \item satisfies the strong version \eqref{eq:CD-inter} of the reduced curvature-dimension inequality
   for $N' = N$ and for all $0 \le s < t < r \le 1$,
  \item has $\mu_t = \rho_tm \in \mathcal{P}^{ac}(X,m)$ for all $t \in [0,1]$ with the density upper-bound
 \begin{equation}\label{eq:densitybound_finite}
  ||\rho_t||_{L^{\infty}(X,m)} \le e^{\sqrt{K^-N}D} \max\{||\rho_0||_{L^{\infty}(X,m)},||\rho_1||_{L^{\infty}(X,m)}\}.
 \end{equation}
 \end{enumerate}
\end{theorem}

\begin{remark}\label{rem:betterclaim}
 As will be seen from the proof, it would suffice in Theorem \ref{thm:main} to assume the inequality \eqref{eq:CDred-def}
 to hold only for $t = 1/2$. Also, we could obtain a slightly better geodesic which also satisfies \eqref{eq:CDred-def}
 for all $N' > N$. See Remark \ref{rem:betterproof} for details how this improvement could be achieved.
\end{remark}

%

The proof of Theorem \ref{thm:main} has two parts. 
The first part is to combine the midpoints to a full geodesics which satisfies the reduced curvature-dimension inequality.
A version of this in nonbranching spaces was proven in \cite{BS2010}. Because we will construct the geodesics by taking
midpoints with minimal entropy, the stronger inequality \eqref{eq:CD-inter} for the critical entropy will follow.

With similar approach as in \cite{R2011b}, we show that the minimizer of the entropy $\sE_N$ among midpoints between 
two given measures has the appropriate density bound.
Here we rely on the fact proven in \cite{R2011b} that the minimizer of a suitable excess mass functional $\mathcal{F}_C$ is zero.

Before proving details for these two parts we will recall some relevant definitions and results on optimal mass transportation in the next section.

\begin{ack}
 Many thanks are due to Nicola Gigli for suggesting this problem and for Karl-Theodor Sturm for useful discussions. 
 This note was written while the author was visiting the Hausdorff Research Institute for Mathematics in Bonn.
 He wishes to thank the Institute for the hospitality.
\end{ack}

\section{Properties of the Wasserstein distance}\label{sec:W}

Let us start with notation which is often used in the theory of optimal mass transportation.
For a comprehensive introduction on the subject see for instance \cite{V2009}.
We denote the set of all Borel probability measures on the space $(X,d)$ by $\mathcal{P}(X)$. For a given Borel measure $m$
on $(X,d)$ we denote by $\mathcal{P}^{ac}(X,m)$ the subset of $\mathcal{P}(X)$ consisting of all the Borel probability
measures that are absolutely continuous with respect to $m$.
Recall that the (possibly infinite) Wasserstein-distance between two Borel probability measures $\mu, \nu \in \mathcal{P}(X)$ is given by
\[
 W_2(\mu, \nu) = \left(\inf\left\{\int_{X\times X} d(x,y)^2\,d\sigma(x,y)\right\}  \right)^{1/2},
\]
where the infimum is taken over all $\sigma \in \mathcal{P}(X \times X)$ with $\mu$ as its first marginal and $\nu$ as the second.

Any geodesic $(\mu_t) \in \Geo(\mathcal{P}(X))$ between two measures $\mu_0, \mu_1 \in \mathcal{P}(X)$
in the space $(\mathcal{P}(X), W_2)$  can be realized as a measure $\pi \in \mathcal{P}(\Geo(X))$ so that $\mu_t = (e_t)_\#\pi$,
where $e_t(\gamma) = \gamma_t$ for any geodesic $\gamma$ and $t \in [0,1]$ and $f_\#\mu$ denotes the push-forward of the measure
$\mu$ under $f$, see for example \cite[Corollary 7.22]{V2009}. Notice that this realization is usually not unique.
We denote by $\GeoOpt(\mu_0, \mu_1)$ the space consisting of all measures $\pi\in \mathcal{P}(\Geo(X))$
for which the mapping $t \mapsto (e_t)_\#\pi$ is a geodesic in $\mathcal{P}(X)$ from $\mu_0 = (e_0)_\#\pi$ to $\mu_1 = (e_1)_\#\pi$.
All geodesics in this note are understood as constant speed geodesics parametrized by the interval $[0,1]$.
We will use the notation $AC([0,1];X)$ for the space of absolutely continuous curves from the interval $[0,1]$ to $X$.

Now we recall some notation and results that were used in \cite{R2011b}.
For any two measures $\mu_0, \mu_1 \in \mathcal{P}(X)$ with $W_2(\mu_0,\mu_1) < \infty$ we define the set
of all the midpoints as
\[
\mathcal{M}(\mu_0,\mu_1) = \left\{\nu \in \mathcal{P}(X) \,:\,  W_2(\mu_0,\nu) =  W_2(\mu_1,\nu) = \frac12 W_2(\mu_0,\mu_1)\right\}. 
\]
We will use the following two basic properties of the set of midpoints. For their proof see \cite{R2011b}. First of all,
for the existence of the minimizers of various functionals we will use the compactness of the set $\mathcal{M}(\mu_0,\mu_1)$
in the Wasserstein space. Notice that $CD^*(K,N)$-spaces are boundedly compact since they are complete and doubling.

\begin{lemma}\label{lma:cmpt}
 Assume that $(X,d)$ is a boundedly compact metric space and that $\mu_0, \mu_1 \in \mathcal{P}(X)$ have bounded support.
 Then the set $\mathcal{M}(\mu_0,\mu_1)$ is compact in $(\mathcal{P}(X),W_2)$.
\end{lemma}

In order to prove the density-bound along the geodesic we will need to redefine a possible bad geodesic in the part where the
density is large. For this we need the following lemma which guarantees that when we redefine part of the geodesic we stay 
inside the set of midpoints.

\begin{lemma}\label{lma:combined}
 Suppose $\mu_0, \mu_1 \in \mathcal{P}(X)$ with $W_2(\mu_0,\mu_1)<\infty$. Then
 for any $\pi \in \GeoOpt(\mu_0,\mu_1)$ and any Borel function $f \colon \Geo(X) \to [0,1]$ with $c = (f\pi)(\Geo(X)) \in (0,1)$ we have
 \[
  (e_{\frac12})_\#\left((1-f)\pi\right) + c\nu \in \mathcal{M}(\mu_0, \mu_1)
 \]
 with every 
 \[
  \nu \in \mathcal{M}\left(\frac1{c} (e_{0})_\#\left(f\pi\right),
      \frac1{c} (e_{1})_\#\left(f\pi\right)\right).
 \]             
\end{lemma}




In \cite{R2011b} the proof for the upper-bound on the densities of measures used the excess mass functional
$\mathcal{F}_C \colon \mathcal{P}(X) \to [0,1]$ defined for all thresholds $C \ge 0$ as
\[
 \mathcal{F}_C(\mu) = ||(\rho-C)^+||_{L^1(X,m)} + \mu^s(X),
\]
where $\mu = \rho m + \mu^s$ with $\mu^s \perp m$, and $a^+ = \max\{0,a\}$.
The existence of the minimizers for this functional follows from the lower semicontinuity, which was
again proved in \cite{R2011b}.

\begin{lemma}\label{lma:lsc}
 Let $(X,d)$ be a bounded metric space with a finite measure $m$.
 Then for any $C \ge 0$ the functional $\mathcal{F}_C$ is lower semicontinuous in $(\mathcal{P}(X),W_2)$.
\end{lemma}

Since this time we will also minimize the entropy $\sE_N$ we will need the lower semicontinuity for it as well.
In addition, the entropy $\sE_N$ may in general attain even the value $-\infty$. However, this is not possible
when the reference measure $m$ is finite. Although our measure $m$ is not necessarily finite, we always consider
measures living in some bounded set and bounded sets have finite $m$-measure. We will skip the proof of the lemma because
the lower semicontinuity is a standard fact and the boundedness away from $-\infty$ is a direct consequence of Jensen's inequality.

\begin{lemma}\label{lma:lscE}
 Let $(X,d)$ be a bounded metric space with a finite measure $m$.
 Then for any $N > 0$ the functional $\sE_N$ is lower semicontinuous in $(\mathcal{P}(X),W_2)$ and attains
 only values on a compact interval.
\end{lemma}


In nonbranching spaces we may decompose a transport to smaller parts which we can then consider separately. Although 
in branching spaces such decomposition is not possible we still have some weaker form of separation.
This will be seen from the next Proposition \ref{prop:separation} and the Lemma \ref{lma:separation} following it.

\begin{proposition}\label{prop:separation}
 Let $\mu_0, \mu_1 \in \mathcal{P}(X)$ with $W_2(\mu_0,\mu_1) < \infty$ and $t_0 \in (0,1)$.
 Suppose that there exist constants $0 \le C_1 \le C_2 < \infty$ and a measure  $\pi \in \GeoOpt(\mu_0,\mu_1)$ with
 \begin{equation}\label{eq:lengthbounds}
  C_1 \le l(\gamma) \le C_2 \qquad\text{for }\pi\text{-a.e. }\gamma \in \Geo(X).
 \end{equation}
 Then the bounds in \eqref{eq:lengthbounds} hold $\tilde\pi$-a.e. for any $\tilde\pi \in \GeoOpt(\mu_0,\mu_1)$ with 
 $(e_{t_0})_\#\tilde\pi = (e_{t_0})_\#\pi$.
\end{proposition}
\begin{proof}
 Take $\tilde\pi \in \GeoOpt(\mu_0,\mu_1)$ with $(e_{t_0})_\#\tilde\pi = (e_{t_0})_\#\pi$.
 Consider the following gluing of the measures $\pi$ and $\tilde\pi$ at $\mu_{t_0} = (e_{t_0})_\#\pi$: 
 Let $\hat\pi \in \mathcal{P}(AC([0,1];X))$ be such that $(e_{t})_\#\pi = (e_{t})_\#\hat\pi$ for $t \in [0,t_0]$
 and $(e_{t})_\#\tilde\pi = (e_{t})_\#\hat\pi$ for $t \in [t_0,1]$. Because we glued together transports via
 the measure $\mu_{t_0}$, the measure $\hat\pi$ is concentrated on the set
 \[
  \left\{\gamma \in  AC([0,1];X) \,:\, \exists w \in [0,1]\text{ with } \gamma^w \in \Geo(X)\right\},
 \]
 where the reparametrization $\gamma^w$ is defined as
 \[
  \gamma^w \colon t \mapsto \begin{cases}
                                              \gamma(\frac{t_0t}{w}), & \text{if }0 \le t \le w\\
                                              \gamma(\frac{t_0-1}{w-1}(t-1)+1), & \text{if }w \le t \le 1.
                                             \end{cases}
 \]

 Let $\hat\pi^r \in \GeoOpt(\mu_0,\mu_1)$ be the reparametrization of $\hat\pi$ given by the set of parameters $w$.
 What needs to be proven is that the reparametrization is trivial, meaning that we can take $w = t_0$ for $\hat\pi$-a.e.
 $\gamma \in  AC([0,1];X)$. This follows by considering the distances between the measure $\mu_{t_0}$ and the endpoints.
 We have
 \begin{align*}
  \int_{\Geo(X)}(w(\gamma)l(\gamma))^2\,d\hat\pi^r(\gamma) & = t_0^2 \int_{\Geo(X)}l(\gamma)^2\,d\hat\pi^r(\gamma),\\
  \int_{\Geo(X)}((1-w(\gamma))l(\gamma))^2\,d\hat\pi^r(\gamma) & = (1-t_0)^2 \int_{\Geo(X)}l(\gamma)^2\,d\hat\pi^r(\gamma),
 \end{align*}
 where we have identified $w(\gamma^w)$ as the $w$ for $\gamma$. These two equalities imply
 \[
  \int_{\Geo(X)}\left(w(\gamma)-t_0\right)\left(w(\gamma)+t_0\right)l(\gamma)^2\,d\hat\pi^r(\gamma) = 0
       = \int_{\Geo(X)}\left(w(\gamma)-t_0\right)l(\gamma)^2\,d\hat\pi^r(\gamma)
 \]
 giving $w = t_0$ for $\hat\pi^r$-almost every $\gamma$.
\end{proof}

In using Proposition \ref{prop:separation} we will need the following consequence of cyclical monotonicity.

\begin{lemma}\label{lma:separation}
 Take $0 \le C_1\le C_2 \le C_3 \le C_4 \le \infty$ and define 
 \[
  A_1 = \{\gamma \in \Geo(X) \,:\, C_1 \le l(\gamma) \le C_2 \}
 \quad \text{and} \quad
  A_2 = \{\gamma \in \Geo(X) \,:\, C_3 < l(\gamma) \le C_4 \}.
 \]
 Then for any $\pi \in \GeoOpt(\mu_0,\mu_1)$ and any $t \in (0,1)$ there exists a set $E \subset \Geo(X)$ with $\pi(E)=0$
 such that
 \[
  \{(\gamma, \hat\gamma) \in (A_1\setminus E) \times (A_2\setminus E)\,:\, \gamma_t = \hat\gamma_t\} = \emptyset.
 \]
\end{lemma}
\begin{proof}
 We know that the optimal plan obtained from $\pi$ is concentrated on a cyclically monotone set, see for instance \cite[Theorem 5.10]{V2009}.
 Suppose that inside this set for some $t \in (0,1)$ we have $\gamma_t = \hat\gamma_t$. Then by cyclical monotonicity we have
 \begin{align*}
 0 \le & ~d(\gamma_0,\hat\gamma_1)^2 + d(\hat\gamma_0,\gamma_1)^2 - d(\gamma_0,\gamma_1)^2 - d(\hat\gamma_0,\hat\gamma_1)^2 \\
   \le & ~d(\gamma_0,\gamma_t)^2 + 2d(\gamma_0,\gamma_t)d(\gamma_t,\hat\gamma_1)+ d(\gamma_t,\hat\gamma_1)^2
         + d(\hat\gamma_0,\gamma_t)^2 + 2d(\hat\gamma_0,\gamma_t)d(\gamma_t,\gamma_1)+ d(\gamma_t,\gamma_1)^2 \\
       & - d(\gamma_0,\gamma_t)^2 - 2d(\gamma_0,\gamma_t)d(\gamma_t,\gamma_1) - d(\gamma_t,\gamma_1)^2 
         - d(\hat\gamma_0,\gamma_t)^2 - 2d(\hat\gamma_0,\gamma_t)d(\gamma_t,\hat\gamma_1) - d(\gamma_t,\hat\gamma_1)^2 \\
    = & ~2(d(\gamma_0,\gamma_t)-d(\hat\gamma_0,\gamma_t))(d(\gamma_t,\hat\gamma_1)-d(\gamma_t,\gamma_1)) \le 0
 \end{align*}
 and so $d(\gamma_0,\gamma_1) = d(\hat\gamma_0,\hat\gamma_1)$.
\end{proof}

 Let us briefly discuss how Proposition \ref{prop:separation} and Lemma \ref{lma:separation} will be used.
 Suppose that we have selected some $\pi \in \GeoOpt(\mu_0,\mu_1)$ and we want to consider separately different
 parts of this measure depending on the lenghts of the curves $\gamma \in \Geo(X)$. Fix any of the measures $\mu_t = (e_t)_\#\pi$
 with $t \in (0,1)$ and consider a geodesic between $\mu_0$ and $\mu_t$ (or $\mu_t$ and $\mu_1$).
 We know by Proposition \ref{prop:separation} that if we redefine part of this geodesic for all curves that have lengths in some
 interval then the new redefined part of the geodesic will live exactly on curves with lengths in this interval.
 By Lemma \ref{lma:separation} this part of the geodesic will be disjoint from the rest of the geodesic for all times except
 possibly for times $0$ and $1$.
 All this means that we can essentially assume the lengths of the curves, and therefore also the distortion coefficients,
 to be constant.

\section{Construction of the good geodesic}

 The construction of the geodesic $(\mu_t)$ in Theorem \ref{thm:main} is done inductively as follows.
 Suppose that for some $n \in \N$ we have defined the measures $\mu_{k2^{-n}}$ for all $0\le k \le 2^n$.
 Then we define for all odd $0 < k < 2^{n+1}$ the measure $\mu_{k2^{-n-1}}$ to be a minimizer of $\sE_N$ in
 $\mathcal{M}(\mu_{(k-1)2^{-n-1}},\mu_{(k+1)2^{-n-1}})$. 
 Such minimizer exists by Lemma \ref{lma:cmpt} and Lemma \ref{lma:lscE}.
 This procedure defines the geodesic for all dyadic times $t$.
 The rest of the geodesic is then given by taking the completion in $(\mathcal{P}(X), W_2)$.
 
 For simplicity we will only consider the case $K \le 0$ in the following computations. The case $K > 0$ follows
 with minor modifications.

\begin{remark}\label{rem:betterproof}
 As mentioned in Remark \ref{rem:betterclaim}, with some more work we could prove a bit stronger result.
 The geodesic $(\mu_t)_{t\in[0,1]}$ in Theorem \ref{thm:main} could also be constructed so that in addition to the
 two estimates stated in the theorem it would also satisfy
 \eqref{eq:CDred-def} for all $t \in [0,1]$ and $N' \ge N$. The main difference is that the minimization of $\sE_N$
 should then be done in a closed subset of $\mathcal{M}(\mu_0,\mu_1)$ where all $\sE_{N'}$, with $N' \ge N$, are bounded from
 above by the value given by the corresponding convexity-inequality.

 However, the main reason why the convexity-inequalities are required for all $N'\ge N$ is to get a definition where $CD^*(K,N_1)$
 implies $CD^*(K,N_2)$ for all $N_2 \ge N_1$. Typically the convexity-inequality is used only for the critical $N$.
 This is why we have decided here to prove only the simpler result stated in Theorem \ref{thm:main}.
\end{remark}

\subsection{Convexity-inequalities for intermediate times}

The geodesic which we constructed from the midpoints satisfies the correct convexity-inequalities because of the
special structure of the coefficients $\sigma_{K,N}^{(t)}(r)$ in the definition of the $CD^*(K,N)$-spaces. One manifestation of this
structure is stated in the following lemma. It is easy to see that the lemma is not true for the coefficients $\tau_{K,N}^{(t)}(r)$
which are used in the definition of $CD(K,N)$-spaces. Coincidentally the construction of a geodesics from the midpoints
in a $CD(K,N)$-space need not produce a geodesic which satisfies the convexity-inequality of the $CD(K,N)$-space.

\begin{lemma}\label{lma:induction}
Let $t_1,t_2, t_3 \in [0,1]$ with $t_1 < t_2$ and let $r \ge 0$. Then
\[
  \sigma_{K,N}^{((1-t_3)t_1+t_3t_2)}(r) =
  \sigma_{K,N}^{(1-t_3)}((t_2-t_1)r)\sigma_{K,N}^{(t_1)}(r) + \sigma_{K,N}^{(t_3)}((t_2-t_1)r)\sigma_{K,N}^{(t_2)}(r)
\]
and
\[
  \sigma_{K,N}^{(1-(1-t_3)t_1-t_3t_2)}(r) =
  \sigma_{K,N}^{(1-t_3)}((t_2-t_1)r)\sigma_{K,N}^{(1-t_1)}(r) + \sigma_{K,N}^{(t_3)}((t_2-t_1)r)\sigma_{K,N}^{(1-t_2)}(r).
\]
\end{lemma}
\begin{proof}
 Let us only prove the first equality as the second one follows by analogous calculations. Let us abbreviate $d = \sqrt{\frac{-K}{N}}r$
 and expand both sides of the claimed equality
\begin{align*}
 \sigma_{K,N}^{((1-t_3)t_1+t_3t_2)}(r) &=  \frac{\sinh(((1-t_3)t_1+t_3t_2)d)}{\sinh(d)}\\
 \sigma_{K,N}^{(1-t_3)}((t_2-t_1)r)&\sigma_{K,N}^{(t_1)}(r) + \sigma_{K,N}^{(t_3)}((t_2-t_1)r)\sigma_{K,N}^{(t_2)}(r) \\
 & = \frac{\sinh((1-t_3)(t_2-t_1)d)\sinh(t_1d)}{\sinh((t_2-t_1)d)\sinh(d)} + \frac{\sinh(t_3(t_2-t_1)d)\sinh(t_2d)}{\sinh((t_2-t_1)d)\sinh(d)}.\\
\end{align*}
 Now the equality follows from basic properties of hyperbolic functions
\begin{align*}
 4\sinh(((1-t_3)&t_1+t_3t_2)d)\sinh((t_2-t_1)d) \\
  =~ & 2\cosh((-t_3t_1+ (1+t_3)t_2)d) - 2\cosh(((2-t_3)t_1 + (t_3-1)t_2)d) \\
  =~ & 2\cosh((t_3t_1+ (1-t_3)t_2)d) - 2\cosh(((2-t_3)t_1 + (t_3-1)t_2)d) \\
 & + 2\cosh((-t_3t_1+ (1+t_3)t_2)d) - 2\cosh((t_3t_1+ (1-t_3)t_2)d) \\
  =~ & 4\sinh((1-t_3)(t_2-t_1)d)\sinh(t_1d) + 4\sinh(t_3(t_2-t_1)d)\sinh(t_2d).
\end{align*}
\end{proof}

Let us now prove that we have the convexity-inequality of the $CD^*(K,N)$-space for our good geodesic.
First we will show that the convexity-inequality holds in the usual form between $0$, $t$ and $1$.
To this aim we first prove the following crude upper-bound for the entropies
\begin{equation}\label{eq:crude}
 \sE_N(\mu_t) \le \sigma_{K,N}^{(1-t)}(D)\sE_N(\mu_0) + \sigma_{K,N}^{(t)}(D)\sE_N(\mu_1),
\end{equation}
which we will eventually only need for showing the continuity of the entropy at the end-points.

 By Proposition \ref{prop:separation} we have for all $n \in \N$ and $0 \le k < 2^n$ that
 $\pi \in \GeoOpt(\mu_{k2^{-n}},\mu_{(k+1)2^{-n}})$ is concentrated on geodesics of length at most $2^{-n}D$.
 Suppose that for some $n \in \N$ the inequality \eqref{eq:crude} is true for all $t=k2^{-n}$ with $0 \le k \le 2^n$.
 Let $0 < k <2^{n+1}$ be odd. Since we assumed that we have a $CD^*(K,N)$-space we have by Lemma \ref{lma:induction}
 \begin{align*}
  \sE_N(\mu_{k2^{-n-1}}) & \le \sigma_{K,N}^{(\frac12)}(2^{-n}D)\sE_N(\mu_{(k-1)2^{-n-1}}) + \sigma_{K,N}^{(\frac12)}(2^{-n}D)\sE_N(\mu_{(k+1)2^{-n-1}}) \\
    & \le \sigma_{K,N}^{(\frac12)}(2^{-n}D)\Big(\sigma_{K,N}^{(1-(k-1)2^{-n-1})}(D)\sE_N(\mu_0) + \sigma_{K,N}^{((k-1)2^{-n-1})}(D)\sE_N(\mu_1)\\
     & + \sigma_{K,N}^{(1-(k+1)2^{-n-1})}(D)\sE_N(\mu_0) + \sigma_{K,N}^{((k+1)2^{-n-1})}(D)\sE_N(\mu_1)\Big) \\
     & \le \sigma_{K,N}^{(1-k2^{-n-1})}(D)\sE_N(\mu_0) + \sigma_{K,N}^{(k2^{-n-1})}(D)\sE_N(\mu_1).
 \end{align*}
 By the lower semicontinuity of $\sE_N$ we then have \eqref{eq:crude} for all $t \in [0,1]$.
 
 Next we prove \eqref{eq:CDred-def} for $N' = N$ and all $t \in [0,1]$. We know that it holds for $t = 1/2$. 
 Suppose that \eqref{eq:CDred-def} holds for $t=k2^{-n}$ with some $n \in \N$ and for all $0 \le k \le 2^n$. Now let
 $n \ge 2$ and $0 < k <2^{n+1}$ be odd. Suppose first that $1 < k <2^{n+1}-1$. Since $\mu_\frac12$ has been fixed,
 by Proposition \ref{prop:separation} and Lemma \ref{lma:separation} we can consider the parts of the transports
 with $d(x_0,x_1) \in [l\epsilon, (l+1)\epsilon[$, $l \in \N$, separately. Therefore, by estimating in these parts similarly as
 for \eqref{eq:crude} and then letting $\epsilon \searrow 0$ the inequality \eqref{eq:CDred-def} follows for $t=k2^{-n-1}$. 

 Let us then consider the case $k \in \{1, 2^{n+1}-1\}$. For $t \in \{0,1\}$ the claim in Lemma \ref{lma:separation} is not valid.
 However, this will not cause real problems since by the crude estimate \eqref{eq:crude} we know that
 \[
  \sE_N(\mu_0) = \lim_{t\searrow 0}\sE_N(\mu_t) \quad \text{and} \quad \sE_N(\mu_1) = \lim_{t\nearrow 1}\sE_N(\mu_t).
 \]
 Therefore it suffices to first do the estimates between $\delta$, $2^{-n-1}$ and $2^{-n}$
 (or between $1-\delta$, $1-2^{-n-1}$ and $1-2^{-n}$ respectively) and then let $\delta \searrow 0$.

 Now that we have obtained the convexity-inequality \eqref{eq:CDred-def} for all $t \in [0,1]$ we prove that also the stronger
 inequality \eqref{eq:CD-inter} holds for our geodesic. Take $0 \le s < r \le 1$. Let $n \in \N$ be the smallest integer for which 
 $s < k2^{-n} < r$ for some $k \in \N$. Then
 \[
  (k-1)2^{-n} \le s < k2^{-n} < r \le (k+1)2^{-n},
 \]
 and so for any measure $\mu \in \mathcal{P}(X)$ with
 \[
  W_2(\mu, \mu_s) = \frac{k2^{-n} - s}{r - s} W_2(\mu_s,\mu_r) \quad \text{and} \quad W_2(\mu, \mu_r) = \frac{r - k2^{-n}}{r - s} W_2(\mu_s,\mu_r)
 \]
 we have $\mu \in \mathcal{M}(\mu_{(k-1)2^{-n}}, \mu_{(k+1)2^{-n}})$. Therefore \eqref{eq:CD-inter} holds with the selected
 $s$ and $r$ and $t = k2^{-n}$.
 
 Then inductively with similar arguments as above \eqref{eq:CD-inter} holds with $t = k2^{-n}$ for all $k,n \in \N$ with
 $s < k2^{-n} < r$. Indeed, assume that $k,n \in \N$ are such that $s < k2^{-n} < r$ and \eqref{eq:CD-inter} holds
 for with all $k' \in \N$ with $s < k'2^{-n+1} < r$. Then, \eqref{eq:CD-inter} holds with $\max\{s,(k-1)2^{-n}\}$,
 $k2^{-n}$ and $\min\{r,(k+1)2^{-n}\}$ and thus by Lemma \ref{lma:induction} with $s$, $k2^{-n}$ and $r$.
 By the lower semicontinuity it finally holds for all $t \in (s,r)$.

\subsection{Upper-bound for the densities}
 
 Our next aim is to prove that for all $t \in [0,1]$ the measure $\mu_t$ is absolutely continuous with respect to $m$
 with a density $\rho_t$ satisfying \eqref{eq:densitybound_finite}.
 As in \cite{R2011b} we first observe how the curvature-dimension condition spreads the support of the measure.

\begin{lemma}\label{lma:spreading}
 Suppose that $(X,d,m)$ is a $CD^*(K,N)$-space with $K \in \R$ and $N \in (1,\infty)$.
 Then for any $\mu_0, \mu_1 \in \mathcal{P}^{ac}(X,m)$ with bounded support and with densities $\rho_0$ and $\rho_1$
 bounded from above there exists $\pi \in \GeoOpt(\mu_0,\mu_1)$ so that
 \begin{equation}\label{eq:bigsupport}
  m(\{x \in X ~:~ \rho_{\frac12}(x) > 0\}) \ge 
  \frac1{e^{\sqrt{K^-N}D/2} \max\{||\rho_0||_{L^{\infty}(X,m)},||\rho_1||_{L^{\infty}(X,m)}\}},
 \end{equation}
 where $(e_\frac12)_\#\pi = \rho_\frac12 m + \mu_\frac12^s$ with $\mu_\frac12^s \perp m$ and 
 $D$ is an upper-bound for the length of $\pi$-almost every $\gamma \in \Geo(X)$.
\end{lemma}
\begin{proof}
   Write
  \[
   M = \max\{||\rho_0||_{L^{\infty}(X,m)},||\rho_1||_{L^{\infty}(X,m)}\}
  \]
  and
  \[
   E = \{x \in X ~:~ \rho_{\frac12}(x) > 0\}.
  \]
  Let $\pi \in \GeoOpt(\mu_0,\mu_1)$ be a measure satisfying \eqref{eq:CDred-def}
  which is concentrated on geodesics with length at most $D$. From \eqref{eq:CDred-def} we get
 \begin{align*}
  \sE_N\left((e_{\frac12})_\#\pi\right) 
   & \le - \iint_{X\times X}\left[\sigma_{K,N'}^{(1/2)}(d(x_0,x_1))\rho_0^{-1/N}(x_0) 
                 + \sigma_{K,N}^{(1/2)}(d(x_0,x_1))\rho_1^{-1/N'}(x_1)\right]\,dq(x_0,x_1)\\
   & \le - e^{-\sqrt{K^-/N}\frac{D}2}M^{-\frac1N},
 \end{align*}
 because we have
 \begin{align*}
   \sigma_{K,N}^{(1/2)}(d(x_0,x_1)) & = \frac{\mathfrak{S}_{K/N}(d(x_0,x_1)/2)}{2\mathfrak{S}_{K/N}(d(x_0,x_1))}
                        = \frac{\sinh(\sqrt{K^-/N}d(x_0,x_1)/2)}{\sinh(\sqrt{K^-/N}d(x_0,x_1))} \\
                       & = \frac{1}{e^{\sqrt{K^-/N}d(x_0,x_1)/2} + e^{-\sqrt{K^-/N}d(x_0,x_1)/2}}
                         \ge \frac12 e^{-\sqrt{K^-/N}\frac{D}2}.
 \end{align*} 
 On the other hand by Jensen's inequality we have
 \[
  \sE_N\left((e_{\frac12})_\#\pi\right) = - \int_E \rho_{\frac12}^{1-\frac1N}\,dm 
                               \ge -m(E)\left(\frac1{m(E)}\int_E \rho_{\frac12}\,dm\right)^{1-\frac1N}
                                \ge - m(E)^{\frac1N}.
 \]
 Combination of these two inequalities gives \eqref{eq:bigsupport}.
\end{proof}

Now continuing from the estimate in Lemma \ref{lma:spreading} exactly as in \cite[Proposition 3.11]{R2011b} 
we obtain the following result.

\begin{proposition}\label{prop:evenspread}
 Assume that $(X,d,m)$ is a $CD^*(K,N)$-space for some $K \in \R$ and $N \in (1,\infty)$
 and that $\mu_0, \mu_1 \in \mathcal{P}^{ac}(X,m)$
 have bounded support and densities $\rho_0$ and $\rho_1$, respectively. Suppose in addition that all
 measures in $\GeoOpt(\mu_0,\mu_1)$ are concentrated on geodesics with length at most $D$. Then we have
 \[
  \min_{\nu \in \mathcal{M}(\mu_0,\mu_1)} \mathcal{F}_C(\nu) = 0
 \]
 for
 \[
  C =  e^{\sqrt{K^-N}D/2}\max\{||\rho_0||_{L^{\infty}(X,m)},||\rho_1||_{L^{\infty}(X,m)}\}.
 \]
\end{proposition}

Next we continue with a result similar to Proposition \ref{prop:evenspread}. In the proof of \cite[Proposition 3.11]{R2011b}
we studied the minimizer of the excess mass functional $\mathcal{F}_C$ using the geodesics satisfying the convexity-inequality
for the entropy $\sE_N$.
This time we will look into the properties of the minimizers of $\sE_N$ and as a tool for spreading the mass we will
use the minimizers of $\mathcal{F}_C$ and the bounds given by the previous Proposition \ref{prop:evenspread}.

\begin{proposition}\label{prop:verygg}
Let $(X,d,m)$ be a $CD^*(K,N)$-space with $K \in \R$ and $N \in (1,\infty)$.
Suppose that $\mu_0, \mu_1 \in \mathcal{P}^{ac}(X,m)$
have bounded support and densities $\rho_0$ and $\rho_1$. Suppose in addition that all
measures in $\GeoOpt(\mu_0,\mu_1)$ are concentrated on geodesics with length at most $D$.
Then for any minimizer $\nu$ of $\sE_N$ in $\mathcal{M}(\mu_0,\mu_1)$ we have 
$\mathcal{F}_C(\nu) = 0$ with
\[
 C =  e^{\sqrt{K^-N}D/2}\max\{||\rho_0||_{L^{\infty}(X,m)},||\rho_1||_{L^{\infty}(X,m)}\}.
\]
\end{proposition}

\begin{proof}
 First of all there exists a minimizer of the entropy $\sE_N$ in $\mathcal{M}(\mu_0,\mu_1)$
 because of Lemma \ref{lma:cmpt} and Lemma \ref{lma:lscE}. Also by Lemma \ref{lma:lscE} the entropy
 at the minimizer is finite. Let $\nu$ be one of the minimizers of $\sE_N$ in $\mathcal{M}(\mu_0,\mu_1)$.

 If $\nu$ would have a singular part with respect to $m$, we could lower the entropy by redefining the 
 part of the geodesic that goes via the singular part to go via a measure satisfying \eqref{eq:CDred-def}.
 This lowering of the entropy would then contradict the minimality of the entropy among the midpoints at $\nu$.
 We can therefore write $\nu = \rho m$. 

 Suppose now, contrary to the claim, that $\mathcal{F}_C(\nu) > 0$.
 Let $\eta > 0$ be such that
 \[
  m(\{x \in X \,:\, \rho(x) > C + \eta\}) \ge m(\{x \in X \,:\, \rho(x) > C + 2\eta\}) > 0.
 \]
 Define
 \[
  C_1 = \frac{4}{\eta}\left(m(\{x \in X \,:\, \rho(x) > C + \eta\})-m(\{x \in X \,:\, \rho(x) > C + 2\eta\})\right).
 \]
 Now for any $\phi \in (0,\frac{\eta}3)$ there exists $\delta \in (\eta,2\eta)$ so that $m(A') < m(A) + C_1\phi$, where
 \[
  A = \{x \in X \,:\, \rho(x) > C + \delta\}\qquad \text{and}\qquad A' = \{x \in X \,:\, \rho(x) \ge C + \delta - 3\phi\}.
 \]

 Take $\pi_1\in \GeoOpt(\nu, \mu_0)$ and $\pi_2\in \GeoOpt(\nu, \mu_1)$, and using Proposition \ref{prop:evenspread} find a measure
 \[
  \tilde\nu = \tilde\rho m \in \mathcal{M}\left(\frac{(e_1)_\# \pi_1|_{\{\gamma_0 \in A\}}}{\nu(A)},\frac{(e_1)_\# \pi_2|_{\{\gamma_0 \in A\}}}{\nu(A)}\right)
 \]
 with
 \begin{equation}\label{eq:estim1}
  \tilde\rho(x) \le \frac{C}{\nu(A)}.
 \end{equation}

Now consider a new measure $\hat\nu = \hat\rho m$ defined as the combination
 \[
  \hat\nu = \nu|_{X \setminus A} + \frac{C + \delta-\phi}{C+ \delta} \nu|_{A} + \frac{\phi}{C+\delta}\nu(A) \tilde\nu.
 \]
By Lemma \ref{lma:combined} we have $\hat\nu \in \mathcal{M}(\mu_0,\mu_1)$.

For $x \in A$ we have the estimates
\begin{align*}
 \hat\rho(x) & \le \frac{C + \delta -\phi}{C+ \delta}\rho(x) + \frac{\phi}{C+\delta}\nu(A)\tilde\rho(x)
  \le \frac{(C + \delta-\phi)\rho(x) + C\phi}{C+\delta} \\
  & = \rho(x) + \frac{(C-\rho(x))\phi}{C + \delta} < \rho(x) - \frac{\delta\phi}{C + \delta}
\end{align*}
and
\[
 \hat\rho(x) \ge \frac{C + \delta - \phi}{C+ \delta}\rho(x) > C + \delta - \phi.
\]

For $x \in A' \setminus A$ we have
\[
 \hat\rho(x) \le \rho(x) + \frac{\phi}{C+\delta}\nu(A)\tilde\rho(x) \le \rho(x) + \frac{C\phi}{C + \delta} < C + \delta + \phi
\]
and for $x \in X \setminus A'$ we get
\[
 \hat\rho(x) \le \rho(x) + \frac{\phi}{C+\delta}\nu(A)\tilde\rho(x) \le C + \delta - 3\phi + \frac{C\phi}{C + \delta} 
  < C + \delta - 2\phi.
\]

Next we need to estimate how much mass is moved from over $A$ to $X \setminus A$.
This is given by
\[
 \kappa_{A} = \int_A (\rho(x) - \hat\rho(x))\,dm  \ge \frac{\delta\phi}{C + \delta}m(A) =: C_2\phi.
\]
Next we estimate the smaller part of this mass which goes on $A' \setminus A$:
\[
 \kappa_{A' \setminus A} = \int_{A'\setminus A} (\hat\rho(x) - \rho(x))\,dm  < C_1\phi^2.
\]
What is left over from $A' \setminus A$ is the larger part which then goes on $X \setminus A'$
\[
 \kappa_{X \setminus A'} = \int_{X\setminus A'} (\hat\rho(x) - \rho(x))\,dm. 
\]

With these we can estimate the change in the entropy when we replace $\nu$ by $\hat\nu$.
\begin{align*}
 \sE_{N}(\hat\mu)-\sE_{N}(\mu)  = &~ -\int_X\hat\rho(x)^{1-\frac1{N}}\,dm + \int_X\rho(x)^{1-\frac1{N}}\,dm \\
  \le &~\int_{A} \hat\rho(x)^{-\frac1{N}}(\rho(x)-\hat\rho(x))\,dm + \int_{X \setminus A'} \hat\rho(x)^{-\frac1{N}}(\rho(x)-\hat\rho(x))\,dm \\
  & + \int_{A' \setminus A} \hat\rho(x)^{-\frac1{N}}(\rho(x)-\hat\rho(x))\,dm\\
  \le &~ \kappa_A\left(C + \delta-\phi\right)^{-\frac1{N}} - \kappa_{X\setminus A'}\left(C + \delta - 2\phi\right)^{-\frac1{N}}
    - \kappa_{A' \setminus A}\left(C + \delta+\phi\right)^{-\frac1{N}}\\
   = &~\kappa_A\left(\left(C + \delta-\phi\right)^{-\frac1{N}} -  \left(C + \delta - 2\phi\right)^{-\frac1{N}}\right) \\
 &    + \kappa_{A' \setminus A}\left(\left(C + \delta - 2\phi\right)^{-\frac1{N}} - \left(C + \delta+\phi\right)^{-\frac1{N}}\right)< 0,
\end{align*}
since there exists a constant $C_3>0$ such that
\begin{align*}
 & \kappa_A\left(\left(\frac{C + \delta-\phi}{C + \delta - 2\phi}\right)^{-\frac1{N}} -  1\right)
   + \kappa_{A' \setminus A}\left(1 - \left(\frac{C + \delta + \phi}{C + \delta - 2\phi}\right)^{-\frac1{N}}\right) \\
 & =\kappa_A\left(\left(1+ \frac{\phi}{C + \delta - 2\phi}\right)^{-\frac1{N}} -  1\right)
   + \kappa_{A' \setminus A}\left(1 - \left(1 + \frac{3\phi}{C + \delta - 2\phi}\right)^{-\frac1{N}}\right) \\
& \le \kappa_A\left(- \frac{\phi}{N(C + \delta - 2\phi)} + \frac{C_3}{N}\phi^2\right)
   + \kappa_{A' \setminus A}\left(\frac{3\phi}{N(C + \delta - 2\phi)}\right) \\
& \le \frac{3C_1\phi^3- C_2\phi^2}{N(C + \delta - 2\phi)} + \frac{C_2C_3}{N}\phi^3 < 0\\
\end{align*}
for small enough $\phi > 0$.
\end{proof}

 Getting the final upper-bound on the density using Proposition \ref{prop:verygg} is easy.
 Suppose that for some $n \in \N$ we have $\mu_{k2^{-n}} = \rho_{k2^{-n}} m$ for all integers $0 \le k \le 2^{n}$ and that for these
 we have
 \begin{equation}\label{eq:indassum}
  ||\rho_{k2^{-n}}||_{L^{\infty}(X,m)} \le 
      \prod_{i=1}^ne^{\sqrt{K^-N}2^{-i}D}\max\{||\rho_0||_{L^{\infty}(X,m)},||\rho_1||_{L^{\infty}(X,m)}\}.
 \end{equation}
 
 Because of the assumption $\diam(\spt \mu_1 \cup \spt \mu_1) = D < \infty$ we have by Proposition \ref{prop:separation}
 that any measure in $\GeoOpt(\Gamma_{k2^{-n}},\Gamma_{(k+1)2^{-n}})$ is concentrated on geodesics with length at most $2^{-n}D$.

 Now for all odd $0 \le k \le 2^{n+1}$ by Proposition \ref{prop:verygg} we have $\mu_{k2^{-n-1}} = \rho_{k2^{-n-1}} m$ with the
 estimate
 \begin{align*}
  ||\rho_{k2^{-n-1}}||_{L^{\infty}(X,m)} 
     & \le e^{\sqrt{K^-N}2^{-n-1}D}\max\{||\rho_{(k-1)2^{-n-1}}||_{L^{\infty}(X,m)},||\rho_{(k+1)2^{-n-1}}||_{L^{\infty}(X,m)}\}\\
     & \le \prod_{i=1}^{n+1}e^{\sqrt{K^-N}2^{-i}D}\max\{||\rho_0||_{L^{\infty}(X,m)},||\rho_1||_{L^{\infty}(X,m)}\}.
 \end{align*} 
 Thus for all dyadic $t\in[0,1]$ the measure $\mu_t$ is absolutely continuous with respect to $m$ and the estimate
 \eqref{eq:densitybound_finite} holds. By the lower semicontinuity of the functionals $\mathcal{F}_C$
 this is then true for all $t \in [0,1]$.


%

%




\begin{thebibliography}{10}
  \bibitem{BS2010} K. Bacher and K.-T. Sturm.
    \emph{Localization and Tensorization Properties of the Curvature-Dimension Condition for Metric Measure Spaces},
     Journal Funct. Anal. \textbf{259} (2010), no. 1, 28--56.
  \bibitem{CS2011} F. Cavalletti and K.-T. Sturm,
    \emph{$CD_\text{loc}(K,N)$ implies $MCP(K,N)$},
    Journal Funct. Anal., to appear.
  \bibitem{DS2011} Q. Deng and K.-T. Sturm,
    \emph{Localization and Tensorization Properties of the Curvature-Dimension Condition for Metric Measure Spaces II},
    Journal Funct. Anal. \textbf{260} (2011), no. 12,  3718--3725.
  \bibitem{LV2009} J. Lott and C. Villani,
     \emph{Ricci curvature for metric-measure spaces via optimal transport},
     Ann. of Math. \textbf{169} (2009), no. 3, 903--991.
   \bibitem{R2011b} T. Rajala,
     \emph{Interpolated measures with bounded density in metric spaces satisfying the curvature-dimension conditions of Sturm},
     preprint.
  \bibitem{R2011} T. Rajala,
     \emph{Local Poincar\'e inequalities from stable curvature conditions on metric spaces},
     Calc. Var. Partial Differential Equations, to appear.
  \bibitem{S2006I} K.-T. Sturm,
     \emph{On the geometry of metric measure spaces. I},
     Acta Math. \textbf{196} (2006), no. 1, 65--131.
  \bibitem{S2006II} K.-T. Sturm,
     \emph{On the geometry of metric measure spaces. II},
     Acta Math. \textbf{196} (2006), no. 1, 133--177.
  \bibitem{V2009} C. Villani,
     \emph{Optimal transport. Old and new},
     vol. 338 of Grundlehren der Mathematischen Wissenschaften, Springer-Verlag, Berlin, 2009.
\end{thebibliography}
\end{document}